\input amstex
\input Amstex-document.sty

\def\p{\partial}

\def\Om{\Omega}

\def\pom{{\p \Om}}
\def\bom{{\overline\Om}}
\def\R{\bold R}

\def\wtt{\widetilde}

\def\Ga{\Gamma}

\def\dist{\text{dist}}
\def\det{\text{det}}

\def\ol{\overline}
\def\lan{\langle}
\def\ran{\rangle}

\def\phi{\varphi}
\def\C{{\Cal C}}

\def\M{\Cal M}
\def\N{\Cal N}


\pageno 221

\def\shorttitle{Affine Maximal Hypersurfaces}
\def\shortauthor{Xu-Jia Wang}

\topmatter
\title\nofrills{\boldHuge Affine Maximal Hypersurfaces}
\endtitle

\author \Large Xu-Jia Wang* \endauthor

\thanks *Centre for Mathematics and Its Applications, The Australian National University, Canberra, ACT 0200,
Australia. E-mail: X.J.Wang\@maths.anu.edu.au
\endthanks

\abstract\nofrills \centerline{\boldnormal Abstract}

\vskip 4.5mm

{\ninepoint This is a brief survey of recent works by Neil Trudinger and myself on the Bernstein problem and
Plateau problem for affine maximal hypersurfaces.

\vskip 4.5mm


}
\endabstract
\endtopmatter

\document

\vskip 2mm

\specialhead \noindent \boldLARGE 1. Introduction \endspecialhead

The concept of affine maximal surface in affine geometry corresponds to that of minimal surface in Euclidean
geometry. The affine Bernstein problem and affine Plateau problem, as proposed in [9,5,7],  are two fundamental
problems for affine maximal surfaces. We shall describe some recent advances, mostly obtained by Neil Trudinger
and myself [21-24], on these two problems.

Given an immersed hypersurface $\M\subset\R^{n+1}$, one defines
the {\it affine metric} (also called the Berwald-Blaschke metric)
by $g=|K|^{-1/(n+2)}I\!I$, where $K$ is the Gauss curvature,
$I\!I$ is the second fundamental form of $\M$. In order that the
metric is positive definite, the hypersurface will always be
assumed to be locally uniformly convex, namely it has positive
principal curvatures. From the affine metric one has the {\it
affine area functional},
$$A(\M) =\int_\M K^{1/(n+2)},  \tag 1.1$$
which can also be written as
$$A(u)=\int_\Om [\det D^2 u]^{1/(n+2)}\tag 1.2$$
if $\M$ is given as the graph of a convex function $u$ over a
domain $\Om\subset\R^n$. The affine metric and affine surface area
are invariant under unimodular affine transformations.

A locally uniformly convex hypersurface is called {\it affine
maximal} if it is stationary for the functional $A$ under interior
convex perturbation. A convex function is called an affine maximal
function if its graph is affine maximal. Traditionally such
hypersurfaces were called affine minimal [1,9]. Calabi suggested
using the terminology affine maximal as the second variation of
the affine area functional is negative [5]. If the hypersurface is
a graph of convex function $u$, then $u$ satisfies the {\it affine
maximal surface equation} (the Euler-Lagrange equation of the
functional $A$),
$$L[u] := U^{ij}w_{ij}=0,\tag 1.3$$
where $[U^{ij}]$ is the cofactor of the Hessian matrix $D^2 u$,
$$w =[\det D^2 u]^{-(n+1)/(n+2)} , \tag 1.4$$
and the subscripts $i, j$ denote partial derivatives with respect
to the variables $x_i, x_j$. Note that for any given $i$ or $j$,
$U^{ij}$, as a vector field in $\Om$, is divergence free. The
equation (1.3) is a nonlinear fourth order partial differential
equation, which can also be written in the short form
$$\Delta_g h=0 , \tag 1.5$$
where $h=(\det D^2 u)^{-1/(n+2)}$, and $\Delta_g$ denotes the
Laplace-Beltrami operator with respect to $g$.

The quantity
$$H_A(\M)=\frac{-1}{n+1} L(u)$$
is called the affine mean curvature of $\M$, and is also invariant
under unimodular affine transformations. In particular it is
invariant if one rotates the coordinates or adds a linear function
to $u$. The affine mean curvature of the unit sphere is $n$.

The affine Bernstein problem concerns the uniqueness of entire
convex solutions to the affine maximal surface equation, and asks
whether an entire convex solution of (1.3) is a quadratic
polynomial. The Chern conjecture [9] asserts this is true in
dimension two. Geometrically, and more generally, it can be stated
as that a Euclidean complete, affine maximal, locally uniformly
convex surface in 3-space must be an elliptic paraboloid. Calabi
proved the assertion assuming in addition that the surface is
affine complete [5], see also [6,7]. A problem raised by Calabi,
called the Calabi conjecture in [19], is whether affine
completeness alone is enough for the Bernstein theorem. The Chern
conjecture was proved true in [21] (see Theorem 3.1 below). The
Calabi conjecture was resolved in [22], as a byproduct of our
fundamental result that affine completeness implies Euclidean
completeness for locally uniformly convex hypersurfaces of
dimensions larger than one (Theorem 3.2). See also [14] for a
different proof of the Calabi conjecture.

The affine Plateau problem deals with the existence and regularity
of affine maximal hypersurfaces with prescribed boundary of which
the normal bundles on the boundary coincide with that of a given
locally uniformly convex hypersurface.  The affine Plateau
problem, which had not been studied before, is more complicated
when compared with the affine Bernstein problem in 3-space. The
first boundary value problem, namely prescribing the solution and
its gradient on the boundary, is a special case of the affine
Plateau problem. We need to impose two boundary conditions as the
affine maximal surface equation is a fourth order equation. We
will formulate the Plateau problem as a variational maximization
problem and prove the existence and regularity of maximizers to
the problem in 3-space [24] (Theorem 5.1). For the existence we
need a uniform cone property of locally convex hypersurfaces,
proved in [23], which also led us to the proof of the conjecture
by Spruck in [20] (Theorem 4.1), concerning the existence of
locally convex hypersurfaces of constant Gauss curvature.

Equation (1.3) can be decomposed as a system of two second order
partial differential equations, one of which is a linearized
Monge-Amp\`ere equation and the other is a Monge-Amp\`ere
equation, see (2.6) and (2.7) below. This formulation enables us
to establish the regularity for equation (1.3) (Theorem 2.1),
using the regularity theory for Monge-Amp\`ere type equations
[2,3]. A crucial assumption in Theorem 2.1 is the strict convexity
of solutions, which is the key issue for both the affine Bernstein
and affine Plateau problems. We succeeded in proving the necessary
convexity estimates only in dimension two.

\vskip 2mm

\specialhead \noindent \boldLARGE 2. A priori estimates \endspecialhead

Instead of the homogeneous equation (1.3), we consider here the non- \linebreak homogeneous (prescribed affine
mean curvature) equation
$$L(u)=f\ \ \ \text{in}\ \ \Om, \tag 2.1$$
where $f$ is a bounded measurable function, and $\Om$ is a
normalized convex domain in $\R^n$. A convex domain is called
normalized if its minimum ellipsoid, that is the ellipsoid with
minimum volume among all ellipsoids containing the domain, is a
unit ball.

Let $u$ be a smooth, locally uniformly convex solution of (2.1)
which vanishes on $\pom$. First we need positive upper and lower
bounds for the determinant $\det D^2 u$. For the upper bound we
have, by constructing appropriate auxiliary function, for any
subdomain $\Om'\subset\subset\Om$, the estimate
$$\sup_{x\in\Om'} \det D^2 u(x)\le C, \tag 2.2$$
where $C$ depends only on $n$, $\dist (\Om', \pom)$, $\sup_\Om
|Du|$, $\sup_\Om f$, and $\sup_\Om |u|$.

For the lower bound we need a key assumption, namely a control on
the strict convexity of solutions, which can be measured by
introducing the {\it modulus of convexity}. Let $v$ be a convex
function in $\Om$. For any $y\in \Om$, $h>0$, denote
$$S_{h,v}(y)=\{x\in\Om\ \big{|}\  v(x)=v(y)+Dv(y)(x-y)+h\} .$$
The modulus of convexity of $v$ is a nonnegative function, defined
by
$$\rho_v (r) = \inf_{y\in\Om} \rho_{v, y}(r),
                                        \ \ \ \ r>0,$$
where
$$\rho_{v, y}(r) = \sup
  \{h\ge 0\ \big{|} \ S_{h,v}(y)\subset B_r(y)\} $$
if  there exists $h\ge 0$ such that $S_{h,v}(y)\subset B_r(y)$,
otherwise we define $\rho_{v, y}(r) = 0$. We have $\rho_v (r)>0$
for all $r>0$ if $v$ is strictly convex.

Let $u$ be a smooth, locally uniformly convex solution of (2.1).
Then we have the following lower bound estimate, for any $\Om'
\subset\subset \Omega$,
$$\inf_{x\in\Om'}\det D^2 u(x)\ge C , \tag 2.3$$
where $C$ depends on $n$, $\dist (\Om', \p\Om)$, $\sup_\Om |Du|$,
$\inf_\Om f$, and $\rho_u$. The proof again can be achieved by
introducing an appropriate auxiliary function.

 From the a priori estimates (2.2) and (2.3) we then have

\proclaim{Theorem 2.1}  Let $u\in C^4(\Om)\cap C^0(\bom)$ be a
locally uniformly convex solution of (2.1). Then for any subdomain
$\Om'\subset\subset\Om$, we have:
\newline
(i) $W^{4, p}$ estimate,
$$\|u\|_{W^{4, p}(\Om')}\le C,  \tag 2.4$$
where $p\in [1, \infty)$, $C$ depends on $n, p, \sup_\Om |f|$,
$\dist(\Om', \pom)$, $\sup_\Om |u|$, and $\rho_u$.
\newline
(ii) Schauder estimate,
$$\|u\|_{C^{4, \alpha}(\Om')}\le C , \tag 2.5$$
where $\alpha\in (0, 1)$, $C$ depends on $n,\alpha,
\|f\|_{C^\alpha(\bom)}$, $\dist(\Om', \pom)$, $\sup_\Om |u|$, and
$\rho_u$.
\endproclaim

Note that the gradient of $u$ is locally controlled by $\rho_u$,
the modulus of convexity of $u$. To prove Theorem 2.1, we write
(2.1) as a second order partial differential system
$$\align
U^{ij}w_{ij} & =f \ \ \ \text{in}\ \ \Om,\tag 2.6\\
\det D^2 u & =w^{-(n+2)/(n+1)}\ \ \ \text{in}\ \ \Om,\tag 2.7\\
\endalign$$
where (2.6) is regarded as a second order elliptic equation for
$w$. By (2.2) and (2.3), and the H\"older continuity of linearized
Monge-Amp\`ere equation [3], we have the interior a priori
H\"older estimate for $w$.  We note that the H\"older continuity
in [3] is proved for the homogeneous equation, but the argument
there can be easily carried over to the non-homogeneous case under
(2.2) and (2.3). By the interior Schauder estimate for the
Monge-Amp\`ere equation [2], we obtain the interior a priori
$C^{2, \alpha}$ estimate for $u$. It follows that (2.6) is a
linear uniformly elliptic equation with H\"older coefficients.
Hence Theorem 2.1 follows.

The control on strict convexity is a key condition in Theorem 2.1.
One cannot expect the strict convexity of solutions when $n\ge 3$.
Indeed, there are convex solutions to the Monge-Amp\`ere equation
$$\det D^2 u=1\tag 2.8$$
which are not strictly convex, and so not smooth [17]. Note that
any non-smooth convex solution of (2.8) can be approximated by
smooth ones, and a smooth solution of (2.8) is obviously a
solution of (2.1), with $f=0$.

An interesting problem is to find appropriate conditions to
estimate the strict convexity of solutions of (2.1). For the
affine Bernstein problem it suffices to prove convexity estimate
for solutions vanishing on the boundary. We succeeded only in
dimension two, see \S 5.

\vskip 2mm

\specialhead \noindent \boldLARGE 3. The affine Bernstein problem \endspecialhead

We say a hypersurface $\M$, immersed in $\R^{n+1}$, is  Euclidean
complete if it is complete under the metric induced from the
standard Euclidean metric.

\proclaim{Theorem 3.1} A Euclidean complete, affine maximal,
locally uniformly convex surface in $\R^3$ is an elliptic
paraboloid.
\endproclaim

Theorem 3.1 extends Jorgens' theorem [11], which asserts that an
entire convex solution of (2.8) in $\R^2$ must be a quadratic
function. Jorgens' theorem also leads to the Bernstein theorem for
minimal surfaces in dimension two [11]. Jorgens' theorem was
extended to higher dimensions by Calabi [4] for $2\le n\le 5$ and
Pogorelov [17] for $n\ge 2$. See also [8]. Observe that the Chern
conjecture follows from Theorem 3.1 immediately.

The proof of Theorem 3.1 uses the affine invariance of equation
(1.3) and the a priori estimates in \S 2. First note that a
Euclidean complete locally uniformly convex hypersurface must be a
graph. Suppose the surface in Theorem 3.1 is the graph of a
nonnegative convex function $u$ with $u(0)=0$. For any constant
$h>1$, let $T_h$ be the linear transformation which normalizes the
section $S^0_{h, u}=\{u<h\}$, and let $v_h(x)=
h^{-1}u(T_h^{-1}(x))$. By the convexity estimate in dimension two,
the modulus of convexity of $v_h$ is independent of $h$. Hence
there is a uniform positive distance from the origin to the
boundary $\p T_h(S^0_{h, u})$. By Theorem 2.1, we infer that the
largest eigenvalue of $T_h$ is controlled by the least one of
$T_h$, which implies that $u$ is defined in the entire $\R^2$. By
Theorem 2.1 again, $D^3 v_h(0)$ is bounded. Hence for any given
$x\in\R^2$,
$$|D^2 u(x)-D^2 u(0)|\le Ch^{-1/2}\to 0$$
as $h\to 0$, namely $D^2 u(x)=D^2 u(0)$.

Note that the dimension two restriction is used only for the
strict convexity estimate. The affine Bernstein problem was
investigated by Calabi in a number of papers [5,6,7]. Using the
result that a nonnegative harmonic function (i.e. $h$ in (1.5))
defined on a complete manifold with nonnegative Ricci curvature
must be a constant, he proved that, among others, the Bernstein
theorem in dimension two, under the additional hypothesis that the
surface is also complete under the affine metric.

Instead of the Euclidean completeness as in the Chern conjecture,
Calabi asks whether affine completeness alone is sufficient for
the Bernstein theorem. This question was recently answered
affirmatively in [22]. See also [14] for a different treatment
based on the result in [16]. In [22] we proved a much stronger
result. That is

\proclaim{Theorem 3.2} An affine complete, locally uniformly
convex hypersurface in $\R^{n+1}$, $n\ge 2$, is also Euclidean
complete.
\endproclaim

The converse of Theorem 3.2 is not true [12], nor is it for $n=1$.
For the proof, which uses the Legendre transform and Lemma 4.1
below, we refer the reader to [22] for details.

\vskip 2mm

\specialhead \noindent \boldLARGE 4. Locally convex hypersurfaces with boundary \endspecialhead

In this section we present some results in [23], which guarantee
the sub-convergence of bounded sequences of locally convex
hypersurfaces with prescribed boundary.

Recall that a hypersurface $\M\subset\R^{n+1}$ (not necessarily
smooth) is called locally convex if it is a locally convex
immersion of a manifold $\N$ and there is a continuous vector
field on the convex side of $\M$, transversal to $\M$ everywhere.
Let $T$ denote the immersion, namely $\M=T(\N)$. For any given
point $x\in\M$, $T^{-1}(x)$ may contains more than one point. To
avoid confusion when referring to a point $x\in\M$ we understand a
pair $(x,p)$ for some point $p\in\N$ such that $x= T(p)$. We say
$\omega_x\subset\M$ is a neighborhood of $x\in\M$ if it is the
image of a neighborhood of $p$ in $\N$. The $r$-neighborhood of
$x$, $\omega_r(x)$, is the connected component of $\M\cap B_r(x)$
containing the point $x$. In [23] we proved the following
fundamental lemma for locally convex hypersurfaces.

\proclaim {Lemma 4.1}    Let $\M$ be a compact, locally convex
hypersurface in $\R^{n+1}$, $n>1$. Suppose the boundary $\p\M$
lies in the hyperplane $\{x_{n+1}=0\}$. Then any connected
component of $\M\cap\{x_{n+1}<0\}$ is convex.
\endproclaim

A locally convex hypersurface $\M$ is called convex if it lies on
the boundary of the convex closure of $\M$ itself. From Lemma 4.1
it follows that a (Euclidean) complete locally convex hypersurface
with at least one strictly convex point is convex, and that a
closed, locally convex hypersurface is convex. Lemma 4.1 also
plays a key role in the proof of Theorem 3.2.

An application we will use here is the uniform cone property for
locally convex hypersurfaces. Let $\C_{x, \xi, r, \alpha}$ denote
the cone
$$\C_{x, \xi, r, \alpha}
  =\{y\in\R^{n+1}\ \big{|}\ |y-x|<r,\ \lan y-x, \xi\ran
   \ge \cos \alpha \, |y-x| \}. $$
We say that $\C_{x, \xi, r, \alpha}$ is an inner contact cone of
$\M$ at $x$ if this cone lies on the concave side of
$\omega_r(x)$. We say $\M$ satisfies the {\it uniform cone
condition} with radius $r$ and aperture $\alpha$ if $\M$ has an
inner contact cone at all points with the same $r$ and $\alpha$.

\proclaim{Lemma 4.2} Let $\M\subset B_R(0)$ be a locally convex
hypersurface with boundary $\p\M$. Suppose $\M$ can be extended to
$\wtt \M$ such that $\p\M$ lies in the interior of $\wtt\M$ and
$\wtt\M-\M$ is locally strictly convex.  Then there exist $r,
\alpha>0$ depending only on $n$, $R$, and the extended part
$\wtt\M-\M$, such that the $r$-neighborhood $\omega_r(x)$ is
convex for any $x\in\M$, and $\M$ satisfies the uniform cone
condition with radius $r$ and aperture $\alpha$.
\endproclaim

In [23] we have shown that if $\p\M$ is smooth and $\M$ is smooth
and locally uniformly convex near $\p\M$, then $\M$ can be
extended to $\wtt\M$ as required in Lemma 4.2. The main point of
Lemma 4.2 is that $r$ and $\alpha$ depend only on $n, R$ and the
extended part $\wtt\M-\M$. Therefore it holds with the same $r$
and $\alpha$ for a family of locally convex hypersurfaces, which
includes all locally uniformly convex hypersurfaces with boundary
$\p\M$, contained in $B_R(0)$, such that its Gauss mapping image
coincides with that of $\M$. For any sequence of locally convex
hypersurfaces in this family, the uniform cone property implies
the sequence converges subsequently and no singularity develops in
the limit hypersurface. This property is the key for the existence
proof of maximizers to the affine Plateau problem. It also plays a
key role for our resolution of the Plateau problem for prescribed
constant Gauss curvature (as conjectured in [20]), see [23]. We
state the result as follows.

\proclaim{Theorem 4.1}  Let $\Ga=(\Ga_1, \cdots, \Ga_n)
\subset\R^{n+1}$ be a smooth disjoint collection of closed
co-dimension two embedded submanifolds. Suppose $\Ga$ bounds a
locally strictly convex hypersurface $S$ with Gauss curvature
$K(S)>K_0>0$. Then $\Ga$ bounds a smooth, locally uniformly convex
hypersurface of Gauss curvature $K_0$.
\endproclaim

If $S$ is a (multi-valued) radial graph over a domain in $S^n$
which does not contain any hemi-spheres, Theorem 4.1 was
established in [10]. Theorem 4.1 has been extended to more general
curvature functions in [18].

\vskip 2mm

\specialhead \noindent \boldLARGE 5. The affine Plateau problem \endspecialhead

First we formulate the affine Plateau problem as a variational
maximization problem. Let $\M_0$ be a compact, connected, locally
uniformly convex hypersurface in $\R^{n+1}$ with smooth boundary
$\Ga=\p \M_0$. Let $S[\M_0]$ denote the set of locally uniformly
convex hypersurfaces $\M$ with boundary $\Ga$ such that the image
of the Gauss mapping of $\M$ coincides with that of $\M_0$. Then
any two hypersurfaces in $S[\M]$ are diffeomorphic. Let $\ol
S[\M_0]$ denote the set of locally convex hypersurfaces which can
be approximated by smooth ones in $S[\M_0]$. Our variational
affine Plateau problem is to find a smooth maximizer to
$$\sup_{\M\in \ol S[\M_0]} A(\M).\tag 5.1$$

To study (5.1) we need to extend the definition of the affine area
functional to non-smooth convex hypersurfaces. Different but
equivalent definitions can be found in [13]. Here we adopt a new
definition introduced in [21, 24], which is also more
straightforward. Observe that the Gauss curvature $K$ can be
extended to a measure on a non-smooth convex hypersurface, and the
measure can be decomposed as the sum of a singular part and a
regular part, $K=K_s+K_r$, where the singular part $K_s$ is a
measure supported on a set of Lebesgue measure zero, and the
regular part $K_r$ can be represented by an integrable function.
We extend the definition of affine area functional (1.1) to
$$A(\M)=\int_\M K_r^{1/(n+2)} .\tag 5.2$$
The affine area functional is upper semi-continuous [13,15]. See
also [21,24] for different proofs.

A necessary condition for the affine Plateau problem is that the
Gauss mapping image of $\M_0$ cannot contain any semi-spheres.
Indeed if $\M$ is affine maximal such that its Gauss mapping image
contains, say, the south hemi-sphere, then the pre-image of the
south hemi-sphere is a graph of a convex function $u$ over a
domain $\Om$ such that $|Du(x)|\to\infty$ as $x\to\pom$. Then
necessarily $\det D^2 u=\infty$ and so $w=0$ on $\pom$. It follows
that $w\equiv 0$ in $\Om$, a contradiction.

\proclaim{Theorem 5.1} Let $\M_0$ be a compact, connected, locally
uniformly convex hypersurface in $\R^3$  with smooth boundary
$\Ga=\p \M_0$. Suppose the image of the Gauss mapping of $\M_0$
does not contain any semi-spheres.  Then there is a smooth
maximizer to (5.1).
\endproclaim

To prove the existence we observe that by the necessary condition,
there exists a positive constant $R$ such that $\M\subset B_R(0)$
for any $\M\in \ol S[\M_0]$. Hence by the uniform cone property,
Lemma 4.2, any maximizing sequence in $\ol S[\M_0]$ is
sub-convergent. The existence of maximizers then follows from the
upper semi-continuity of the affine area functional. Note that the
existence is true for all dimensions.

To prove the regularity we need to show that
\newline
(i) $\M$ can be approximated by smooth affine maximal surfaces;
and
\newline
(ii) $\M$ is strictly convex.
\newline
The purpose of (i) is such that the a priori estimate in Section 2
is applicable. Note that (i) also implies the Bernstein Theorem
3.1 holds for non-smooth affine maximal surfaces.

By the penalty  method we proved (i) for all dimensions, using the
following classical solvability of the second boundary value
problem for the affine maximal surface equation.

\proclaim{Theorem 5.2} Consider the problem
$$\align
  L(u) & =f(x, u)\ \ \ \text{in}\ \ \Om,\tag 5.3\\
     u & =\phi\ \ \ \ \ \text{on}\ \ \p\Om,\\
     w & =\psi\ \ \ \ \ \text{on}\ \ \p\Om,\\
\endalign $$
where $w$ is given in (1.4), $\Om$ is a uniformly convex domain
with $C^{4, \alpha}$ boundary, $0<\alpha<1$, $f$ is H\"older
continuous, non-decreasing in $u$, \  $\phi, \psi\in C^{4,
\alpha}(\bom)$, and $\psi$ is positive. Then there is a unique
uniformly convex solution $u\in C^{4, \alpha}(\bom)$ to the above
problem.
\endproclaim

To prove Theorem 5.2 we first prove that $u$ satisfies (2.2) and
(2.3),  and that $w$ is Lipschitz continuous on the boundary,
namely $|w(x)-w(y)|\le C|x-y|$ for any $x\in\Om$, $y\in\pom$.
Theorem 5.2 is then reduced to the boundary $C^{2, \alpha}$
estimate for the Monge-Amp\`ere equation. The proof for the
boundary $C^{2, \alpha}$ estimate involves a delicate iteration
scheme. We refer to [24] for details.

The interior $C^{2, \alpha}$ estimate for the Monge-Amp\`ere
equation was proved by Caffarelli [2], using a perturbation
argument. The boundary $C^{2, \alpha}$ estimate, which also uses a
similar perturbation argument, contains substantial new
difficulty, as that the sections
$$S^0_{h,v}(y)=\{x\in\Om\ \big{|}\  v(x)<v(y)+Dv(y)(x-y)+h\}$$
can be normalized for the interior estimate but not for the
boundary estimate. We need to prove that $S^0_{h,v}(y)$ has a good
shape for sufficiently small $h>0$ and $y\in\pom$.

Finally we would like to mention our idea of proving the strict
convexity, namely (ii) above. Note that for both the affine
Bernstein and affine Plateau problem, the dimension two assumption
is only used for the proof of the strict convexity.  To prove the
strict convexity we suppose to the contrary that $\M$ contains a
line segment. Let $P$ be a tangent plane of $\M$ which contains
the line segment. Then the contact set $F$, namely the connected
set of $P\cap\M$ containing the line segment, is a convex set. If
$F$ has an extreme point which is an interior point of $\M$, by
rescaling and choosing appropriate coordinates we obtain a
sequence of affine maximal functions which converges to a convex
function $v$, such that $v(0)=0$ and $v(x)>0$ for $x\ne 0$ in an
appropriate coordinate system, and $v$ is not $C^1$ at the origin
$0$. In dimension two this means $\det D^2 v$ is unbounded near
$0$, which is in contradiction with the estimate (2.2).

If all extreme points of $F$ are boundary points of $\M$, we use
the Legendre transform to get a new convex function which is a
maximizer of a variational problem similar to (5.1), and satisfies
the properties as $v$ above, which also leads to a contradiction.

\vskip 2mm

\specialhead \noindent \boldLARGE 6. Remarks \endspecialhead

We proved the affine Bernstein problem in dimension two. In high
dimensions ($n\ge 10$) a counter-example was given in [21], where
we proved that the function
$$u(x)=(|x'|^9+x_{10}^2)^{1/2}\tag 6.1$$
is affine maximal, where $x'=(x_1, \cdots, x_9)$.

The function $u$ in (6.1) contains a singular point, namely the
origin. The graph of $u$ is indeed an {\it affine cone}, that is
all the level sets $S_{h, u}=\{u=h\}$ are affine self-similar, in
the sense that there is an affine transformation $T_h$ such that
$T_h(S_{h,u})=S_{1,u}$. The above counter-example shows that there
is an affine cone in dimensions $n\ge 10$ which is affine maximal
but is not an elliptic paraboloid. We have not been successful in
finding smooth counter-examples. Little is known for dimensions
$3\le n\le 9$.

For the affine Plateau problem, an interesting problem is whether
the maximizer satisfies the boundary conditions. If $\M_0$ is the
graph of a smooth, uniformly convex function $\phi$, defined in a
bounded domain $\Om\subset\R^n$,  the Plateau problem becomes the
first boundary value problem, that is equation (1.3) subject to
the boundary conditions:
$$\align
   u & = \phi \ \ \text{on}\ \ \pom,\tag 6.2\\
  Du & = D\phi\ \ \text{on}\ \ \pom. \tag 6.3\\
  \endalign$$
In this case we also proved the uniqueness of maximizers of (5.1).
Obviously the maximizer $u$ satisfies (6.2). Whether $u$ satisfies
(6.3) is still unknown. Recall that the Dirichlet problem of the
minimal surface equation is solvable for any smooth boundary
values if and only if the boundary is mean convex. Therefore an
additional condition may be necessary in order that (6.3) is
fulfilled.

\vskip 2mm

\specialhead \noindent \boldLARGE References \endspecialhead


\item {[1]}   W. Blaschke,
               Vorlesungen \"uber Differential geometrie, Berlin, 1923.

\item {[2]}  L.A. Caffarelli,
               Interior $W^{2,p}$ estimates for solutions
               of Monge-Amp\`ere equations,
               Ann. Math., 131 (1990), 135--150.

\item {[3]}  L.A. Caffarelli and C.E. Guti\'errez,
               Properties of the solutions of the linearized
               Monge-Amp\`ere equations,
               Amer. J. Math., 119(1997), 423--465.

\item {[4]} E. Calabi,
               Improper affine hypersurfaces of convex type and a
               generalization of a theorem by K. J\"orgens,
               Michigan Math. J., 5(1958), 105--126.

\item {[5]} E. Calabi,
               Hypersurfaces with maximal affinely invariant area,
               Amer. J. Math. 104(1982), 91--126.

\item {[6]} E. Calabi,
             Convex affine maximal surfaces,
             Results in Math., 13(1988), 199--223.

\item {[7]} E. Calabi,
               Affine differential geometry and holomorphic curves,
               Lecture Notes Math. 1422(1990), 15--21.

\item {[8]} S.Y. Cheng and S.T. Yau,
               Complete affine hypersurfaces, I.
               The completeness of affine metrics,
               Comm. Pure Appl. Math., 39(1986), 839--866.

\item {[9]}  S.S. Chern,
               Affine minimal hypersurfaces,
               in {\it minimal submanifolds and \linebreak geodesics},
               Proc. Japan-United States Sem., Tokyo, 1977, 17--30.

\item {[10]}  B. Guan and J. Spruck,
               Boundary value problems on $S\sp n$ for surfaces
               of constant Gauss curvature.
               Ann. of Math., 138(1993), 601--624.

\item {[11]}   K. Jorgens,
               \"Uber die L\"osungen der Differentialgleichung
               $rt-s\sp 2=1$, Math. Ann. 127(1954), 130--134.

\item {[12]}  K. Nomizu and T. Sasaki,
               Affine differential geometry,
               Cambridge, 1994.

\item {[13]}   K. Leichtweiss,
               Affine geometry of convex bodies,
               Johann Ambrosius Barth Verlag, Heidelberg, 1998.

\item {[14]}  A.M. Li and F. Jia,
               The Calabi conjecture on affine maximal surfaces,
               Result. Math., 40(2001), 265--272.

\item {[15]}   E. Lutwak,
               Extended affine surface area,
               Adv. Math., 85(1991), 39--68.

\item {[16]}  A. Martinez and F. Milan,
               On the affine Bernstein problem,
               Geom. Dedicata, 37(1991), 295--302.

\item {[17]}   A.V. Pogorelov,
               The muitidimensional Minkowski problems,
               J. Wiley, New York, 1978.

\item {[18]}  W.M. Sheng, J. Urbas, and X.-J. Wang,
               Interior curvature bounds for a class of curvature
               equations, preprint.

\item {[19]}   U. Simon,
               Affine differential geometry,
               in {\it Handbook of differential geometry},
               North-Holland, Amsterdam, 2000, 905--961.

\item {[20]}   J. Spruck,
               Fully nonlinear elliptic equations and applications
               in geometry, Proc. International Congress Math.,
               Birkh\"auser, Basel, 1995,   1145--1152.

\item {[21]} N.S. Trudinger and X.-J. Wang,
               The Bernstein problem for affine maximal hypersurfaces,
               Invent. Math., 140 (2000), 399--422.

\item {[22]} N.S. Trudinger and X.-J. Wang,
               Affine complete locally convex hypersurfaces,
               Invent. Math., to appear.

\item {[23]} N.S. Trudinger and X.-J. Wang,
               On locally convex hypersurfaces with boundary,
               J. Reine Angew. Math., to appear.

\item {[24]} N.S. Trudinger and X.-J. Wang,
               The Plateau problem for affine maximal hypersurfaces,
               Preprint.

\enddocument